\documentclass[11pt]{article}
\usepackage{amssymb}
   \csname @addtoreset\endcsname{equation}{section}
\newtheorem{theorem}{Theorem}%[section]
\newtheorem{lemma}{Lemma}%[section]
\newtheorem{corollary}{Corollary}%[section]
\newtheorem{definition}{Definition}%[section]
\newtheorem{remark}{Remark}%[section]
\def\e{\varepsilon}

\def\defi{\stackrel{{\scriptscriptstyle \Delta}}{=}}

\def\a{\alpha}
\def\d{\delta}
\def\o{\omega}
\def\O{\Omega}

\def\F{{\cal F}}
\def\w{\widehat}
\def\Ind{{\,\rm Ind\,}}
\def\Ind{{\mathbb{I}}}

\def\Re{{\rm Re\,}}
\def\Im{{\rm Im\,}}

\def\R{{\bf R}}

\def\L{L}

\def\g{\gamma}
\def\C{{\bf C}}

\def\X{{\cal X}}

\def\M{{\cal M}}

\def\L{{\cal L}}

\newcommand{\be}{\begin{equation}}
\newcommand{\ee}{\end{equation}}
\newcommand{\bd}{\begin{displaymath}}
\newcommand{\ed}{\end{displaymath}}
\newcommand{\ba}{\begin{array}{ll}}
\newcommand{\ea}{\end{array}}
\newcommand{\baa}{\begin{eqnarray}}
\newcommand{\eaa}{\end{eqnarray}}
\newcommand{\baaa}{\begin{eqnarray*}}
\newcommand{\eaaa}{\end{eqnarray*}}
\font\sm=cmr10
\def\a{\alpha}

\def\K{\cal K}
\date{ }
\title{
Predictability of band-limited, high-frequency, and mixed
processes in the presence of ideal low-pass filters}
\author{
Nikolai Dokuchaev\\ {\sm  Department of Mathematics, Trent
University, Ontario, Canada}}
 
\begin{document}
 \vspace{-0.5cm}      \maketitle
\begin{abstract} Pathwise predictability
 of continuous time processes is studied
in deterministic setting. We discuss  uniform prediction  in some
weak sense with respect to certain classes of inputs. More
precisely, we study possibility of approximation of  convolution
integrals over future time by integrals over past time.
 We found that all band-limited processes
are predictable in this sense, as well as high-frequency processes
with zero energy at low frequencies. It follows that a process of
mixed type still can be predicted if an ideal low-pass  filter
exists for this process.
\\    {\bf Key words}: prediction, ideal low-pass filters, band-limited processes,
Hardy spaces, causal estimators.
\\ AMS 2000 classification : 60G25,  93E10, 42B30.
\\
PACS 2008 numbers: 02.30.Mv, %    Approximations and expansions
02.30.Nw,  %  Fourier analysis
02.30.Yy, %    Control theory
07.05.Mh,  %  Neural networks, fuzzy logic, artificial intelligence
07.05.Kf %Data analysis: algorithms and implementation; data management
\end{abstract}
\section{Introduction}  We study pathwise predictability of
continuous time processes in deterministic setting. It is well
known that certain restrictions on frequency distribution can
ensure additional opportunities for prediction and interpolation
of the processes. The classical result is
Nyquist-Shannon-Kotelnikov interpolation theorem for the low-band
processes.  There are related predictability results for low-band
processes (see, e.g., Wainstein and Zubakov (1962), Beutler
(1966), Brown(1969),
 Slepian (1978), Knab
(1981), Papoulis (1985), Marvasti (1986), Vaidyanathan (1987),
Lyman {\it et al} (2000, 2001)).
\par
In the present paper, we study a special kind of weak
predictability such that convolution integrals over future can be
approximated by convolution integrals over past times representing
historical observations. We found some cases when this
approximation can be made uniformly over a wide class of input
processes. We found that all band-limited processes are
predictable in this sense. Similar result is obtained for
high-frequency processes.
 For the
processes of mixed type, we found that the similar predictability
can be achieved when the model allows a low pass filter  that acts
as an ideal low-pass filter for this process. These results can be
a useful addition to the existing theory of band-limited
processes. The novelty is that we consider predictability of both
high frequent and band-limited processes in a weak sense uniformly
over classes of input processes. In addition, we suggest a new
type of predictor. Its kernel is given explicitly in the frequency
domain.
\section{Definitions}
Let  $\Ind$ denote the indicator function, $\R^+\defi[0,+\infty)$,
$\C^+\defi\{z\in\C:\ \Re z>  0\}$, $i=\sqrt{-1}$.
\par
For complex valued functions $x\in L_1(\R)$ or $x\in L_2(\R)$, we
denote by $X=\F x$ the function defined on $i\R$ as the Fourier
transform of $x$; $$(\F x)(i\o)= \int_{-\infty}^{\infty}e^{-i\o
t}x(t)dt,\quad \o\in\R.$$ If $x\in L_2(\R)$, then $X$ is defined
as an element of $L_2(\R)$.
\par
For $v(\cdot)\in L_2(\R)$ such that $v(t)=0$ for $t<0$, we
denote by $\L v$  the Laplace transform \baa\label{Up} V(p)=(\L
v)(p)\defi\int_{0}^{\infty}e^{-p t}v(t)dt, \quad p\in\C^+. \eaa
\par
Let  $H^r$ be the Hardy space of holomorphic on $\C^+$ functions
$h(p)$ with finite norm
$\|h\|_{H^r}=\sup_{s>0}\|h(s+i\o)\|_{L_r(\R)}$, $r\in[1,+\infty]$
(see, e.g., Duren (1970)).
\par
Let $\O>0$ be given.
\begin{definition}
Let $\K$ be the class of functions $k:\R\to\R$ such that $k(t)=0$
for $t>0$ and such that $K=\F k$ is \be
K(i\o)=\frac{d(i\o)}{\d(i\o)}, \label{kda}\ee where $d(\cdot)$ and
$\d(\cdot)$ are polynomials such that ${\rm deg\,} d < {\rm deg\,}
\d$, and if $\d(p)=0$ for $p\in\C$ then $\Re p>0$, $|\Im p|<\O$.
\end{definition}
Note that the class $\K$ is quite wide: it consists of linear
combinations of functions $q(t)e^{\lambda t}\Ind_{\{t\le 0\}}$,
where $\lambda\in\C$, $\Re\lambda>0$, $|\Im \lambda|<\O$, and
where $q(t)$ is a polynomial.
\begin{definition}
Let $\w\K$ be the class of functions $\w k:\R\to\R$ such that $k
(t)=0$ for $t<0$ and such that $K(\cdot)=\L\w k\in H^2\cap
H^\infty$.
\end{definition}
 We are going to study linear predictors in the form $\w y(t)=\int_{-\infty}^t\w
k(t-s)x(s)ds$ for the processes  $y(t)=\int_{t}^{+\infty}
k(t-s)x(s)ds$, where $k\in\K$ and $\w k\in\w K$. The predictors
use historical values of currently observable process $x(\cdot)$.
\par
\begin{definition}\label{def1}
Corrections:

On Endnote 1: it should be  "Let $V$ be as defined above"

On Endnote 2: it should be  "are such as required in statements
(ii)--(iii)".

On Endnote 3: it should be  "Let $V$ and $\widehat K$ be as
defined above"

Additional correction: P.6, Line 25:  please put "$\omega\in
D_{\epsilon}$" instead of "$\omega\in D$".

P.7, Line 9: "Academic Press" instead of "Academic"

P.7, Line 18: "Fourier analysis, and uncertainty" instead of
"Fourier analysis and uncertainty"

Let  $\X=\{x(\cdot)\}$ be a class of functions $x:\R\to \C$. Let
$r\in[1,+\infty]$.
\begin{itemize}
\item[(i)]
 We say that the class $\X$ is  $L_r$-predictable in
the weak sense  if, for any $k(\cdot)\in\K$,  there exists a
sequence $\{\w k_m(\cdot)\}_{m=1}^{+\infty}=\{\w
k_m(\cdot,\X,k)\}_{m=1}^{+\infty}\subset \w\K$ such that $$ \|y-\w
y_m\|_{L_r(\R)}\to 0\quad \hbox{as}\quad m\to+\infty\quad\forall
x\in\X, $$ where \baaa y(t)\defi \int_t^{+\infty}k(t-s)x(s)ds,\qquad
\w y_m(t)\defi \int^t_{-\infty}\w k_m(t-s)x(s)ds.\label{predict}
\eaaa
\item[(ii)] Let the set $\F(\X)\defi \{X=\F x,\quad x\in \X\}$  be provided with a norm $\|\cdot\|$.
 We say that the class $\X$ is  $L_r$-predictable in
the weak sense  uniformly with respect to the norm $\|\cdot\|$, if,
for any $k(\cdot)\in\K$ and $\e>0$, there exists $\w k(\cdot) =\w k
(\cdot,\X,k,\|\cdot\|,\e)\in \w\K$ such that $$ \|y- \w
y\,\|_{L_r(\R)}\le \e\|X\|\quad \forall x\in\X,\quad X=\F x. $$ Here
$y(\cdot)$ is the same as above, $\w y(t)\defi \int^t_{-\infty}\w
k(t-s)x(s)ds.$
\end{itemize}
\end{definition}
 We call functions $\w k(\cdot)$ in Definition \ref{def1} predictors or
predicting kernels.
\section{The main result}
\par
Let $\O>0$ be the same as in the definition of $\K$, and let \baaa
\X_L\defi\{x(\cdot)\in L_2(\R):\
X(\o)=0\quad\hbox{if}\quad|\o|>\O,\quad X=\F x\},\\
\X_H\defi\{x(\cdot)\in L_2(\R):\
X(\o)=0\quad\hbox{if}\quad|\o|<\O,\quad X=\F x\}. \eaaa
 In particular, $\X_L$ is a
class of band-limited processes, and $\X_H$ is a class of
high-frequency processes.
\subsection{Predictability of band-limited and high-frequency
processes from $L_2$}
\begin{theorem}\label{ThM}\begin{itemize}\item[(i)] The classes $\X_L$ and
$\X_H$ are $L_2$-predictable in the weak sense. \item[(ii)] The
classes $\X_L$ and $\X_H$  are $L_\infty$-predictable in the weak
sense  uniformly with respect to the norm $\|\cdot\|_{L_2(\R)}$.
\item[(iii)] For any $q>2$, the classes $\X_L$ and $\X_H$  are
$L_2$-predictable  in the weak sense uniformly with respect to the
norm $\|\cdot\|_{L_q(\R)}$.
\end{itemize}
\end{theorem}
\begin{remark} Since the constant $\O$ is the same for the
classes $\K$, $\X_L$, $\X_H$,  the set of $k(\cdot)\in\K$ such
that the corresponding processes $y(\cdot)$ can be predicted is
restricted for $x(\cdot )\in \X_H$. On the other hand, these
restrictions are absent for band-limited processes $x(\cdot)\in
\X_L$, since they are automatically included to all similar
classes with larger $\O$, i.e., the constant $\O$ in the
definition of $\X_L$ can always be
 increased.
\end{remark}
\par  The question arises how to find the predicting kernels. In
the proof of Theorem \ref{ThM}, a possible choice of the kernels
is given explicitly via  Fourier transforms.
\subsection{Predictability for some bounded processes} Let $C(\R)$ be the
Banach space of all bounded and continuous functions $f:\R\to\C$,
and
 let $C(\R)^*$ be the dual space for $C(\R)$, i.e., it is the space of all linear
continuous functionals $\xi:C(\R)\to\C$ (see, e.g., Yosida
(1980)).
\par
Let  $\M_{\infty}$ be the class of all processes
 $x(t):\R\to\C$ such that there exists a function $X_c\in L_1(\R)$, a
 sequence  $\{\o_k\}_{k=1}^{+\infty}\subset \R$, and
 a sequence $\{c_k\}_{k=1}^{+\infty}\subset\C$ such that
 $\sum_{k=1}^{+\infty}|c_k|<+\infty$ and
 \baaa
 x(t)=\frac{1}{2\pi}\sum_{k=1}^{\infty}c_ke^{i\o_kt}+\frac{1}{2\pi}\int_{-\infty}^{+\infty}e^{i\o
 t}X_c(\o)d\o.
 \eaaa
 Clearly, any set $X\defi
 \left(\{\o_k\}_{k=1}^{+\infty},\{c_k\}_{k=1}^{+\infty},X_c\right)$
 with the required properties is uniquely defined by the process
 $x\in\M_{\infty}$, and
 can be associated with an unique element of $C(\R)^*$ such that
\baaa \langle
f,X\rangle=\sum_{k=1}^{\infty}c_kf(\o_k)+\int_{-\infty}^{+\infty}f(\o)X_c(\o)d\o
 \quad\forall f\in C(\R). \eaaa
In particular,
 $x(t)=\langle
\frac{1}{2\pi} e^{it\cdot},X\rangle$ for all $t$. We will denote
this relationship as $X=\F x$, using the same notation as for the
Fourier transform, and we extend Definition \ref{def1} on this
case  (it is a frequency representation, but not a Fourier
transform anymore). As required in Definition \ref{def1}, we
provide the set $\{X\}$ of these sets $X$ with the norm
$\|\cdot\|_{C(\R)^*}$.
\par
If $x\in\M_{\infty}$, then $ |x(t)|\le
(2\pi)^{-1}\|e^{it\cdot}\|_{C(\R)}\|X\|_{C(\R)^*}$. Hence all
functions from $\M_{\infty}$ are bounded on $\R$.
\par
 Let $\e\in (0,\O)$ be given.
Let \baaa &&\M_L\defi\Bigl\{x\in \M_{\infty}:
|\o_k|\le\O-\e\,\,(\forall k), \quad {\rm supp\,} X_c\subseteq
[-\O+\e,\O-\e]\Bigr\},
\\ &&\M_H\defi\Bigl\{x\in \M_{\infty}: |\o_k|\ge\O+\e\,\,(\forall k), \quad {\rm
supp\,} X_c\subseteq (-\infty,-\O-\e]\cup [\O+\e,+\infty)\}. \eaaa
 $\M_L$ is
a class of band-limited processes, and $\M_H$ is a class of
high-frequency processes.
\begin{theorem}\label{ThM1} The classes $\M_L$ and
$\M_H$ are  $L_\infty$-predictable in the weak sense uniformly
with respect to the norm $\|\cdot\|_{C(\R)^*}$.
\end{theorem}
\section{On a model with ideal low pass-pass filter}
\begin{corollary}\label{cor1}
Assume a model with  a process $x(\cdot)$ such that an observer is
able to decompose it
 as $x(t)=x_L(t)+x_H(t)$, where   $x_L(\cdot)\in\X_L\cup \M_L$ and
$x_H(\cdot)\in\X_H\cup \M_H$. Then this observer would be able to
predict (approximately, in the sense of weak presdictability) the
values of $y(t)= \int_t^{+\infty}k(t-s)x(s)ds$ for $k(\cdot)\in\K$
by predicting the processes $y_L(t)=\int_t^{+\infty}k(t-s)x_L(s)ds$
and $y_H(t)=\int_t^{+\infty}k(t-s)x_H(s)ds$ separately. More
precisely, the process $\w y(t)\defi \w y_L(t)+\w y_H(t)$ is the
prediction of $ y(t)$, where $y_L(t)=\int^t_{-\infty}\w
k_L(t-s)x_L(s)ds$ and $y_H(t)=\int^t_{-\infty}\w k_H(t-s)x_H(s)ds$,
and where $\w k_L(\cdot)$ and $\w k_H(\cdot)$ are predicting kernels
which existence for the processes $x_L(\cdot)$ and $x_H(\cdot)$  is
established above.
\end{corollary}
\par
Let $\chi_L(\o)\defi \Ind_{\{|\o|\le \O\}}$ and $\chi_H(\o)\defi
1-\chi_L(\o)=\Ind_{\{|\o|> \O\}}$, where $\o\in\R$.
\par
 The assumptions of Corollary \ref{cor1} mean
that there are a low-pass filter  and a high-pass filter with the
transfer functions $\chi_L$ and $\chi_H$ respectively, with
$x(\cdot)$ as the input, i.e., that the values $x_L(s)$ and
$x_H(s)$ for $s\le t$ are available at time $t$, where \baaa
x_L(\cdot)\defi \F^{-1}X_L,\quad X_L(\o)\defi \chi_
L(\o)X(\o),\quad x_H(\cdot)\defi \F^{-1}X_H,\quad X_H(\o)\defi
\chi_H(\o)X(\o),\eaaa and where  $X\defi \F x$.
 It
follows that the predictability in the weak sense described in
Definition \ref{def1} is possible  for any process $x(\cdot)$ that
can be decomposed without error on a band limited process and a
high-frequency process, i.e., when there is a low-pass filters
which behaves as an ideal filter for this process. (Since
$x_H(t)=x(t)-x_L(t)$, existence of the law pass filter implies
existence of the high pass filter). On the other hand, Corollary
\ref{cor1} implies that the existence of ideal low-pass filters is
impossible for general processes, since they cannot be predictable
in the sense of Definition \ref{def1}.
\par
Clearly, processes $x(\cdot)\in\X_L\cup \X_H\cup \M_L\cup \M_H$ are
automatically covered by Corollary \ref{cor1}, i.e., the existence
of the filters is not required for this case. For instance, we have
immediately that $x_L(\cdot)=x(\cdot)$ and $x_H(\cdot)\equiv 0$ for
band-limited processes.
\section{Proofs}
Let  $k(\cdot)\in\K$ and $K(i\o)=\F k$. Let (\ref{kda}) holds with
$\d(p)=\prod_{m=1}^n\d_m(p)$, where $\d_m(p)\defi p-a_m+b_mi$, and
where $a_m,b_m\in\R$, $p\in\C$. By the assumptions on $\K$, we
have that $a_m>0$ and $|b_m|<\O$.
\par
It suffices to present a set of predicting kernels $\w k$ with
desired properties. We will use a version of the construction
introduced in Dokuchaev (1996) for an optimal control problem.  This
construction is very straightforward and does not use the advanced
theory of $H^p$-spaces.
\par For
$\g\in\R$, set \baaa &&\a_m=\frac{\O^2-b_m^2}{a_m},\quad
V_m(p)\defi
1-\exp\left(\g\frac{p-a_m+b_mi}{p+\a_m-b_mi}\right),\quad
V(p)\defi \prod_{m=1}^nV_m(p),\\\quad
  &&\w K(i\o)\defi V(i\o)K(i\o).
  \eaaa
\begin{lemma}\label{lemmaV}
\begin{itemize}
\item[(i)] $V(p)\in H^2\cap H^{\infty}$ and $\w K(p)\defi K(p)V(p)\in H^2\cap
H^{\infty}$;
\item[(ii)]  If $\g>0$ and $\o\in[-\O,\O]$, then  $|V(i\o)|\le 1$.
If $\g<0$, and if $\o\in\R$, $|\o|\ge \O$, then $|V(i\o)|\le 1$.
\item[(iii)] If $\o\in(-\O,\O)$, then $V(i\o)\to 1$ as
$\g\to+\infty$. If $\o\in \R$ and  $|\o|>\O$, then $V(i\o)\to 1$
as $\g\to -\infty$.
\item[(iv)]
For any $\e>0$, $V(i\o)\to 1$ as $\g\to+\infty$ uniformly in
$\o\in[-\O+\e,\O-\e]$ as $\g\to+\infty$, and  $V(i\o)\to 1$  as
$\g\to -\infty$ uniformly in $\o\in \R$ such that $|\o|\ge \O+\e$.
\end{itemize}
\end{lemma}
\par
{\it Proof of Lemma \ref{lemmaV}}. Clearly, $V_m(p)\in
H^{\infty}$, and $\d_m(p)^{-1}V_m(p)\in H^2\cap H^{\infty}$, since
the pole of $\d_m(p)^{-1}$ is being compensated by multiplying on
$V_m(p)$. It follows that $K(p)V(p)\in H^2\cap H^{\infty}$. Then
statement (i) follows.
\par Further, for $\o\in\R$,
 \baaa
&&\frac{i\o-a_m+b_mi}{i\o+\a_m-b_mi}=\frac{(-a_m+i\o+ib_m)(\a_m-i\o+b_mi)}{(\o-b_m)^2+a_m^2}\\
&&=\frac{-a_m\a_m+(\o+b_m)(\o-b_m)}{(\o-b_m)^2+\a_m^2}
+i\frac{-a_m(\o+b_m)+\a_m(\o+b_m)}{(\o-b_m)^2+\a_m^2}\,. \eaaa
 Then
 \baaa
\Re\frac{i\o-a_m+b_mi}{i\o+\a_m-b_mi}
=\frac{-a_m\a_m+\o^2-b_m^2}{(\o-b_m)^2+\a_m^2}=\frac{\o^2-\O^2}{(\o-b_m)^2+\a_m^2}\,.
\label{re} \eaaa
 Then statements (ii)-(iv)  follow.
 This completes the proof of Lemma \ref{lemmaV}. $\Box$
\par
{\it Proof of Theorem \ref{ThM}}.  For  $x(\cdot)\in L_2(\R)$, let
 $X\defi \F x$,
 $Y\defi \F y=K(i\o)X(\o)$.
 Let $V$ be as defined above. Set $\w Y(\o)\defi \w K(i\o)X(\o)=V(i\o)Y(\o)$.
\par
 Let us consider the cases of $\X_L$ and $\X_H$ simultaneously.
For the case of the class $\X_L$, consider $\g>0$ and assume that
$\g>0$ and $\g\to +\infty$. Set $D\defi [-\O,\O]$ for this case.
For the case of the class $\X_H$, consider $\g<0$ and assume that
$\g<0$ and $\g\to -\infty$. Set $D\defi
(-\infty,-\O]\cup[\O,+\infty)$ for this case.
\par
Let $x(\cdot)\in\X_L$ or $x(\cdot)\in\X_H$. In  both  cases, Lemma
\ref{lemmaV} gives that $|V(i\o)|\le 1$ for all $\o\in D$. If
$\g\to +\infty$ or $\g\to -\infty$ respectively  for $\X_L$ or
$\X_H$ cases, then
  $V(i\o)\to 1$  for a.e. $\o\in D$, i.e., for a.e. $\o$ such that
  $X(\o)\neq 0$.
\par
Let us prove (i).  Since $K(i\o)\in L_{\infty}(\R)$ and $X\in
L_2(\R)$, we have that $Y(\o)=K(i\o)X(\o)\in L_{2}(\R)$ and $\w Y\in
L_{2}(\R)$.  By Lemma \ref{lemmaV}, it follows that  \baa\w Y(\o)\to
Y(\o)\quad\hbox{for a.e.}\quad \o\in \R, \label{YY}\eaa
 as $\g\to +\infty$
or $\g\to -\infty$ respectively for $\X_L$ or $\X_H$ cases. We
have that $X\in L_2(\R)$, $K(i\o)\in L_2(\R)\cap L_{\infty}(\R)$
and \baa&&|\w K(i\o)-K(i\o)|\le |V(i\o)-1||K(i\o)|\le 2
|K(i\o)|,\quad \o\in D,\label{d1}\\ &&|\w Y(\o)-Y(\o)|\le
2|Y(\o)|=2|K(i\o)||X(\o)|,\quad \o\in D.\label{d2} \eaa
 By (\ref{YY}),(\ref{d2}), and by Lebesque Dominance Theorem, it follows that
\be\|\w Y-Y\|_{L_2(\R)}\to 0,\quad\hbox{i.e.,}\quad\|\w
y-y\|_{L_2(\R)}\to 0 \label{1s}\ee
 as $\g\to +\infty$
or $\g\to -\infty$ respectively for $\X_L$ or $\X_H$ cases, where
$\w y=\F^{-1}\w Y$.
\par
Let us prove (ii)-(iii).
 Take $d=1$ for (ii) and take $d=2$ for (iii). If
$X\in L_q(\R)$ for $q>d$, then H\"older inequality gives
 \be\|\w Y-Y\|_{L_d(\R)}\le
\|\w K(i\o)-K(i\o)\|_{L_\mu(D)}\|X\|_{L_q(D)}, \label{4s}\ee where
$\mu$ is such that $1/\mu+1/q=1/d$.
 By (\ref{d1}) and by Lebesque Dominance Theorem again, it follows that
\be\label{2s}\|\w K(i\o)-K(i\o)\|_{L_\mu(D)}\to 0\quad \forall \mu
\in[1,+\infty) \ee  as $\g\to +\infty$ or $\g\to -\infty$
respectively for $\X_L$ or $\X_H$ cases.
 By (\ref{4s})-(\ref{2s}), it follows  that the predicting kernels $\w
k(\cdot)=\w k(\cdot,\g)=\F^{-1}\w K(i\o)$ are such as required in
statements (ii)--(iii). This completes the proof of Theorem
\ref{ThM}. $\Box$
\par {\it Proof of Theorem \ref{ThM1}}.  For
$x(\cdot)\in \M_{\infty}$ such that  $X=(\{\o_k\}_{k=1}^{+\infty},
\{c_k\}_{k=1}^{+\infty},X_c)$, we have that the corresponding set
$Y=\F y$ is  $Y=\left(\{\o_k\}_{k=1}^{+\infty},
\{K(i\o_k)c_k\}_{k=1}^{+\infty},K(i\o)X_c(\o)\right)$. Similarly
to $X$, it can be considered as an element of $C(\R)^*$ such that
$y(t)=\langle \frac{1}{2\pi}e^{it\cdot},Y\rangle$.
 Let $V$ and $\w K$ be as defined above. Set
  \baaa
  \w Y\defi \left(\{\o_k\}_{k=1}^{+\infty},
\{\w K(i\o_k)c_k\}_{k=1}^{+\infty},\w K(i\o)X_c(\o)\right).
 \eaaa
It can be seen as an element of $C(R)^*$, and $\w
y(t)=\int_{-\infty}^t\w k(t-s)x(s)ds=\langle
\frac{1}{2\pi}e^{it\cdot},\w Y\rangle$, where the kernel is
defined via inverse Fourier transform $\w k(\cdot)=\F^{-1}\w
K(i\o)$.
\par
 We consider the cases of $\M_L$ and $\M_H$ simultaneously.
For the case of the class $\M_L$, we consider $\g>0$ and $\g\to
+\infty$. Set $D_\e\defi[-\O+\e,\O-\e]$ for this case. For the
case of the class $\M_H$, we consider $\g<0$ and  $\g\to -\infty$.
Set $D_\e\defi(-\infty,-\O-\e]\cup[\O+\e,+\infty)$ for this case.
\par
Let $x(\cdot)\in\M_L$ or $x(\cdot)\in\M_H$. In  both  cases, Lemma
\ref{lemmaV} gives that $|V(i\o)|\le 1$ for all $\o\in D_\e$. If
$\g\to +\infty$ or $\g\to -\infty$ respectively  for $\M_L$ or
$\M_H$ cases, then
  $V(i\o)\to 1$  uniformly in $\o\in D_\e$.  Hence $\|\w
K-K\|_{L_{\infty}(D_\e)}\to 0$ as $\g\to +\infty$ or $\g\to
-\infty$, for the cases of $\M_L$ and $\M_H$, respectively. If
$x\in\M_L$ or $x\in\M_H$, then \baaa |\langle f,X\rangle|\le
\max_{t\in D_\e} |f(t)|\,\|X\|_{C(\R)^*}\quad \forall f\in C(\R),
\quad X=\F x. \eaaa Hence  \baaa |\w y(t)-y(t)|=
\left|\left\langle \frac{1}{2\pi}e^{it\cdot},\w
Y-Y\right\rangle\right| =\left|\left\langle
\frac{1}{2\pi}e^{it\cdot}(\w K-K),X\right\rangle\right|\le
\frac{1}{2\pi}\|\w K-K\|_{L_{\infty}(D_\e)}\|X\|_{C(\R)^*} \eaaa
for all $t\in\R$.
 Then the proof of
Theorem \ref{ThM1} follows. $\Box$
\par Corollary \ref{cor1} follows
immediately from Theorem \ref{ThM}.
\begin{remark}  Formally, the
corresponding predictors require the past values of $x(s)$ for all
$s\in(-\infty,t]$, but it is not too restrictive, since
$\int_{-\infty}^t\w k(t-s)x(s)ds$ can be approximated by
$\int_{-M}^t\w k(t-s)x(s)ds$ for large enough $M>0$. In addition,
the corresponding transfer functions can be approximated by
rational fraction polynomials, and more general kernels $k$ can be
approximated by kernels from $\K$.
\end{remark}
\begin{remark}
 The system for the suggested
predictors is stable, since the corresponding transfer functions
have poles in the domain $\{\Re z<0\}$ only.  However, the
suggested predictors are not robust. For instance, if the
predictor is designed for the class $\X_L$ and it is applied for a
process $x(\cdot)\notin\X_L$ with small non-zero energy at the
frequencies outside $[-\O,\O]$, then the error generated by the
presence of this energy is increasing if $\g\to\infty$.
\end{remark}

\end{document}